\pdfoutput=1
\documentclass[a4paper,11pt]{amsart}

\usepackage[margin=1in]{geometry}
\usepackage[utf8]{inputenc}
\usepackage[T1]{fontenc}
\usepackage{lmodern,microtype}
\usepackage{amssymb,mathtools}
\usepackage{enumitem}
\usepackage{tikz}
\usepackage[unicode,bookmarks=false]{hyperref}

\hypersetup{colorlinks=true,linkcolor=black,citecolor=black,filecolor=black,urlcolor=black}

\makeatletter
\let\citationorig\citation
\def\citation#1{\citationorig{#1}\@for\@tempa:=#1\do{\@ifundefined{cit@\@tempa}{\global\@namedef{cit@\@tempa}{}}{}}}
\let\bibitemorig\bibitem
\def\bibitem#1{\@ifundefined{cit@#1}{\typeout{LaTeX Warning: Unused bibitem `#1'}}{}\bibitemorig{#1}}
\makeatother

\renewenvironment{itemize}{\begin{itemorig}[label=\textbullet, noitemsep, topsep=3pt plus 3pt, labelsep=.6em, labelindent=.2em, leftmargin=*]}{\end{itemorig}}

\def\calF{\mathcal{F}}

\let\leq\leqslant

\hypersetup{
  pdftitle={Realization of shift graphs as disjointness graphs of 1-intersecting curves in the plane},
  pdfauthor={Torsten Mütze, Bartosz Walczak, Veit Wiechert},
}

\title{Realization of shift graphs as disjointness graphs of~1-intersecting curves in the plane}
\author{Torsten Mütze\and Bartosz Walczak\and Veit Wiechert}

\begin{document}

\maketitle
\thispagestyle{empty}

The \emph{disjointness graph} of a family $\calF$ of curves in the plane is the graph with vertex set $\calF$ and with edges between the pairs of disjoint curves.
In the \emph{shift graph} $H_m$, the vertices are the ordered pairs $(i,j)$ satisfying $1\leq i<j\leq m$, and two such pairs $(i,j)$ and $(k,\ell)$ form an edge if and only if $j=k$ or $\ell=i$.
The graphs $H_m$ are triangle-free and have unbounded chromatic number \cite[Theorem~6]{EH64} (in fact, it is well known that $\chi(H_m)=\lceil\log_2m\rceil$).
We strengthen a result of Pach, Tardos, and Tóth \cite[Theorem~5]{PTT17} showing that $H_m$ can be realized as the disjointness graph of a family $\calF_m$ of polygonal curves in the plane with the following properties:
(a) every curve in $\calF_m$ is made of $3$ straight-line segments, and
(b) any two curves in $\calF_m$ intersect in at most one point.

Let $\prec$ denote the lexicographic order on the vertices of $H_m$: $(i,j)\prec(k,\ell)$ when $i<k$ or $i=k$ and $j<\ell$.
First, we construct a realization in which the curve $c(i,j)$ representing the vertex $(i,j)$ of $H_m$ consists of $4$ straight-line segments as illustrated, where\hfill
\begin{tikzpicture}[scale=0.5]
  \useasboundingbox (0,3.9) rectangle (6,4);
  \draw (0,1)--node[midway,left]{$A$}(0,0)--node[midway,above]{$B$}(6,0)--node[midway,above right]{$C$}(2,4)--node[midway,left]{$D$}(2,2);
\end{tikzpicture}
\begin{itemize}
\item the $A$-parts occur in the order $\prec$ from right to left,
\item the $B$-parts occur in the order $\prec$ from bottom to top,
\item the $C$-parts occur in the order $\prec$ from top-right to bottom-left,
\item the $D$-parts occur in the order $\prec$ from left to right,
\item the $A$-part of $c(i,j)$ intersects the $B$-parts of $c(k,\ell)$ up to $c(j-1,m)$,
\item the $D$-part of $c(i,j)$ intersects the $B$-parts of $c(k,\ell)$ down to $c(j+1,j+2)$ (if $j\leq m-2$).\hfill
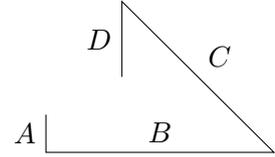
\begin{tikzpicture}[scale=1.2]
  \useasboundingbox (4.1,4.5) rectangle (4.2,4.6);
  \draw (-0.2,0.5)--(-0.2,0)--(4.2,0)--(0,4.2)--(0,0.9);
  \draw (-0.4,0.9)--(-0.4,0.2)--(3.7,0.2)--(0.2,3.7)--(0.2,1.1);
  \draw (-0.6,1.1)--(-0.6,0.4)--(3.2,0.4)--(0.4,3.2)--(0.4,1.1);
  \draw (-0.8,0.9)--(-0.8,0.6)--(2.7,0.6)--(0.6,2.7)--(0.6,1.1);
  \draw (-1,1.1)--(-1,0.8)--(2.2,0.8)--(0.8,2.2)--(0.8,1.1);
  \draw (-1.2,1.1)--(-1.2,1)--(1.7,1)--(1,1.7)--(1,1.1);
  \node[left] at (-0.2,0) {$\scriptscriptstyle(1,2)$};
  \node[left] at (-0.4,0.2) {$\scriptscriptstyle(1,3)$};
  \node[left] at (-0.6,0.4) {$\scriptscriptstyle(1,4)$};
  \node[left] at (-0.8,0.6) {$\scriptscriptstyle(2,3)$};
  \node[left] at (-1,0.8) {$\scriptscriptstyle(2,4)$};
  \node[left] at (-1.2,1) {$\scriptscriptstyle(3,4)$};
  \node[above right] at (2.2,2.2) {\large $H_4$};
\end{tikzpicture}
\end{itemize}

\parshape 1
  0pt .65\textwidth
\noindent
If $(i,j)\prec(k,\ell)$, then $c(i,j)$ and $c(k,\ell)$ intersect if and only if $j\neq k$, as required.
By a standard stretching argument (see e.g.\ \cite[Theorem~3]{MP92}), the $AB$-parts of the curves can be replaced by single straight-line segments with the same intersection pattern.

\parshape 9
  0pt .65\textwidth
  0pt .65\textwidth
  0pt .65\textwidth
  0pt .65\textwidth
  0pt .52\textwidth
  0pt .52\textwidth
  0pt .52\textwidth
  0pt .52\textwidth
  0pt \textwidth
The curves in the non-stretched realization can be arranged to touch (but not cross) a common straight line.
By contrast, triangle-free intersection graphs of curves touching (but not crossing) a common straight line have bounded chromatic number \cite{McG00}.
An analogous contrast but in the opposite way holds for graphs realized by straight-line segments: while there exist triangle-free segment intersection graphs with arbitrarily large chromatic number \cite{PKK+14}, the proof technique of Larman et~al.\ \cite[Theorem~1]{LMPT94} shows that disjointness graphs of straight-line segments (or, more generally, $x$-monotone curves) in the plane satisfy $\chi=O(\omega^4)$.

\end{document}